\documentclass[11pt]{article}
\usepackage{epsfig}
\def\be{\begin{equation}}
\def\ee{\end{equation}}
\def\bea{\begin{eqnarray}}
\def\eea{\end{eqnarray}}
\def\bes{\begin{eqnarray*}}
\def\ees{\end{eqnarray*}}

\def\nn{\nonumber}
\def\<{\langle}
\def\>{\rangle}
\def\lb{\label}
\def\bs{\setminus}

\def\R{{\bf R}}
\def\C{{\bf C}}
\def\Z{{\bf Z}}
\def\N{{\bf N}}
\def\U{{\bf U}}
\def\Q{{\bf Q}}

\def\ga{{\gamma}}

\def\th{{\theta}}
\def\om{{\omega}}
\def\Om{{\Omega}}
\def\ep{{\epsilon}}
\def\lm{{\lambda}}

\def\sg{{\sigma}}

\def\Sg{{\Sigma}}

\def\H{{\cal H}}

\def\P{{\cal P}}
\def\J{{\cal J}}

\def\im{{\rm im}}

\def\Sp{{\rm Sp}}

\def\dm{{\rm \diamond}}

\def\ol#1{\overline{#1}}  %overline in math mode
\def\td#1{\tilde{#1}}
\def\hb{\vrule height0.18cm width0.14cm $\,$}

\def\ol#1{\overline{#1}}  %overline in math mode
\def\td#1{\tilde{#1}}

\title{Stability of closed characteristics on compact hypersurfaces \\in $\R^{2n}$ under pinching condition}

\author{Wei Wang\thanks{Partially supported by LMAM in Peking University
in China and China Postdoctoral Science Foundation No.20070420264.
E-mail: alexanderweiwang@yahoo.com.cn, wangwei@math.pku.edu.cn  }\\
School of Mathematical Science \\ Peking University, Beijing 100871 \\
PEOPLES REPUBLIC OF CHINA \\ }

\date{Jun. 22th, 2008 }

\begin{document}

\maketitle

\begin{abstract}
{\it In this article, let $\Sigma\subset\R^{2n}$ be a compact convex
hypersurface which is $(r, R)$-pinched with
$\frac{R}{r}<\sqrt{\frac{3}{2}}$.
Then  $\Sg$ carries at least two
strictly elliptic closed characteristics; moreover,
$\Sg$ carries at least $2[\frac{n+2}{4}]$
non-hyperbolic closed characteristics.}
\end{abstract}



{\bf Running title}: Stability of closed characteristics

\renewcommand{\theequation}{\thesection.\arabic{equation}}
\renewcommand{\thefigure}{\thesection.\arabic{figure}}

\setcounter{equation}{0}%\setcounter{figure}{0}
\section{Introduction and main results}%{Section 1}

In this article, let $\Sigma$ be a fixed $C^3$ compact convex hypersurface
in $\R^{2n}$, i.e., $\Sigma$ is the boundary of a compact and strictly
convex region $U$ in $\R^{2n}$. We denote the set of all such hypersurfaces
by $\H(2n)$. Without loss of generality, we suppose $U$ contains the origin.
We consider closed characteristics $(\tau,y)$ on $\Sigma$, which are
solutions of the following problem
\be \left\{\matrix{\dot{y}=JN_{\Sigma}(y), \cr
               y(\tau)=y(0), \cr }\right. \lb{1.1}\ee
where $J=\left(\matrix{0 &-I_n\cr
                I_n  & 0\cr}\right)$,
$I_n$ is the identity matrix in $\R^n$, $\tau>0$ and $N_\Sigma(y)$ is
the outward normal vector of $\Sigma$ at $y$ normalized by the
condition $N_{\Sigma}(y)\cdot y=1$. Here $a\cdot b$ denotes the
standard inner product of $a, b\in\R^{2n}$. A closed characteristic
$(\tau,\, y)$ is {\it prime} if $\tau$ is the minimal period of $y$.
Two closed characteristics $(\tau,\, y)$ and $(\sigma, z)$ are {\it
geometrically distinct}  if $y(\R)\not= z(\R)$. We denote by
$\J(\Sg)$ and $\widetilde{\J}(\Sg)$ the set of all closed
characteristics $(\tau,\, y)$ on $\Sg$ with $\tau$ being the minimal
period of $y$ and the set of all geometrically distinct ones
respectively. Note that
$\J(\Sg)=\{\theta\cdot y\,|\, \theta\in S^1,\;y\; is\; prime\}$,
while $\widetilde{\J}(\Sg)=\J(\Sg)/S^1$, where the natural $S^1$-action
is defined by $\theta\cdot y(t)=y(t+\tau\theta),\;\;\forall \theta\in S^1,\,t\in\R$.

Let $j: \R^{2n}\rightarrow\R$ be the gauge function of $\Sigma$, i.e.,
$j(\lambda x)=\lambda$ for $x\in\Sigma$ and $\lambda\ge0$, then
$j\in C^3(\R^{2n}\setminus\{0\}, \R)\cap C^0(\R^{2n}, \R)$
and $\Sigma=j^{-1}(1)$. Fix a constant $\alpha\in(1,\,+\infty)$ and
define the Hamiltonian function
$H_\alpha :\R^{2n}\rightarrow [0,\,+\infty)$ by
\be H_\alpha(x)=j(x)^\alpha,\qquad \forall x\in\R^{2n}.\lb{1.2}\ee
Then
$H_\alpha\in C^3(\R^{2n}\setminus\{0\}, \R)\cap C^1(\R^{2n}, \R)$
is convex and $\Sigma=H_\alpha^{-1}(1)$.
It is well known that problem (\ref{1.1}) is equivalent to
the following given energy problem of Hamiltonian system
\be
\left\{\matrix{\dot{y}(t)=JH_\alpha^\prime(y(t)),
             &&\quad H_\alpha(y(t))=1,\qquad \forall t\in\R. \cr
     y(\tau)=y(0). && \cr }\right. \lb{1.3}\ee
Denote by $\mathcal{J}(\Sigma, \,\alpha)$ the set of all solutions
$(\tau,\, y)$ of (\ref{1.3}) where $\tau$ is the minimal period of
$y$ and by $\widetilde{\mathcal{J}}(\Sigma, \,\alpha)$ the set of
all geometrically distinct solutions  of (\ref{1.3}). As above,
$\widetilde{\mathcal{J}}(\Sigma, \,\alpha)$ is obtained from
$\mathcal{J}(\Sigma, \,\alpha)$ by dividing the natural
$S^1$-action. Note that elements in $\mathcal{J}(\Sigma)$ and
$\mathcal{J}(\Sigma, \,\alpha)$ are one to one correspondent to each
other, similarly for $\widetilde{\J}(\Sg)$ and
$\widetilde{\mathcal{J}}(\Sigma, \,\alpha)$.

Let $(\tau,\, y)\in\mathcal{J}(\Sigma, \,\alpha)$. The fundamental
solution $\gamma_y : [0,\,\tau]\rightarrow \Sp(2n)$ with $\gamma_y(0)=I_{2n}$
of the linearized Hamiltonian system
\be \dot w(t)=JH_\alpha^{\prime\prime}(y(t))w(t),\qquad \forall t\in\R,\lb{1.4}\ee
is called the {\it associate symplectic path} of $(\tau,\, y)$.
The eigenvalues of $\gamma_y(\tau)$ are called {\it Floquet multipliers}
of $(\tau,\, y)$. By Proposition 1.6.13 of \cite{Eke3}, the Floquet multipliers
with their multiplicities of $(\tau,\, y)\in\mathcal{J}(\Sigma)$ do not depend on
the particular choice of the Hamiltonian function in (\ref{1.3}).
For any $M\in \Sp(2n)$, we define the {\it elliptic height } $e(M)$ of
$M$ to be the total algebraic multiplicity of all eigenvalues of $M$ on the
unit circle $\U=\{z\in\C|\; |z|=1\}$ in the complex plane $\C$.
Since $M$ is symplectic, $e(M)$ is even and $0\le e(M)\le 2n$.
As usual $(\tau,\, y)\in\J(\Sg)$ is {\it elliptic} if
$e(\gamma_y(\tau))=2n$. It is {\it strictly elliptic} if
all the eigenvalues $\neq1$ are Krein-definite.
It is {\it non-degenerate} if $1$ is a double
Floquet multiplier of it. It is  {\it hyperbolic} if $1$ is a
double Floquet multiplier of it and $e(\gamma_y(\tau))=2$.
It is well known that these concepts are independent of the choice of $\alpha$.

As in Definition 5.1.7 of \cite{Eke3}, a $C^3$ hypersurface
$\Sigma$ bounding a compact convex region $U$, containing
$0$ in its interior is {\it $(r, R)$-pinched}, with
$0<r\le R$ if
\be |y|^2R^{-2}\le\frac{1}{2}H_2^{\prime\prime}(x)y\cdot y
\le |y|^2r^{-2},\qquad \forall x\in\Sg,\quad \forall y\in\R^{2n}.
\lb{1.5}\ee

For the existence and multiplicity of geometrically distinct closed
characteristics on convex compact hypersurfaces in $\R^{2n}$ we refer to
\cite{Rab1}, \cite{Wei1}, \cite{EkL1}, \cite{EkH1}, \cite{Szu1},
\cite{Vit1}, \cite{HWZ}, \cite{LoZ1}, \cite{LLZ}, \cite{WHL}, and references therein.

On the stability problem, in \cite{Eke2} of Ekeland in 1986 and \cite{Lon2}
of Long in 1998, for any $\Sg\in\H(2n)$ the existence of at least one
non-hyperbolic closed characteristic on $\Sg$ was proved provided
$^\#\td{\J}(\Sg)<+\infty$. Ekeland proved also in \cite{Eke2} the existence
of at least one strictly elliptic closed characteristic on $\Sg$ provided $\Sg\in\H(2n)$
is $\sqrt{2}$-pinched. In \cite{DDE1} of 1992, Dell'Antonio, D'Onofrio and
Ekeland proved the existence of at least one elliptic closed characteristic
on $\Sg$ provided $\Sg=-\Sg$. In \cite{Lon4} of 2000,
Long proved that $\Sg\in\mathcal{\H}(4)$ and $\,^{\#}\td{\J}(\Sg)=2$ imply that both of
the closed characteristics must be elliptic. In \cite{LoZ1} of 2002, Long and
Zhu further proved when $^\#\td{\J}(\Sg)<+\infty$, there exists at least one
elliptic closed characteristic and there are at least $[\frac{n}{2}]$ geometrically
distinct closed characteristics on $\Sg$ possessing irrational mean indices,
which are then non-hyperbolic. In \cite{LoW1}, Long and the author
proved that there exist at least two non-hyperbolic closed characteristics on
$\Sg\in\H(6)$ when $^\#\td{\J}(\Sg)<+\infty$.
In \cite{W1}, the author proved that on every $\Sg\in\H(6)$ satisfying
$^\#\td{\J}(\Sg)<+\infty$, there exist at least two
closed characteristics possessing irrational mean indices and if
$^\#\td{\J}(\Sigma)=3$, then there exist at least two elliptic closed
characteristics. In \cite{W2}, the author studies stability
of closed characteristics on symmetric hypersurfaces.

Motivated by these results, we prove
the following results in this article:

{\bf Theorem 1.1.} {\it let $\Sigma\subset\R^{2n}$ be a compact convex
hypersurface which is $(r, R)$-pinched with
$\frac{R}{r}<\sqrt{\frac{3}{2}}$. Then
$\Sg$ carries at least two geometrically distinct
strictly elliptic closed characteristics. }

{\bf Theorem 1.2.} {\it let $\Sigma\subset\R^{2n}$ be a compact convex
hypersurface which is $(r, R)$-pinched with
$\frac{R}{r}<\sqrt{\frac{3}{2}}$. Then
$\Sg$ carries at least $2[\frac{n+2}{4}]$
geometrically distinct non-hyperbolic closed
characteristics.}

The proof of these theorems is motivated by the methods in \cite{BTZ1},
\cite{Eke3} and \cite{LoZ1} by using the index iteration theory
and comparison theorems on indices as in the study of closed geodesics.

In this article, let $\N$, $\N_0$, $\Z$, $\Q$, $\R$ and $\C$ denote
the sets of natural integers, non-negative integers, integers, rational
numbers, real numbers and complex numbers respectively.
Denote by $a\cdot b$ and $|a|$ the standard inner product and norm in
$\R^{2n}$. Denote by $\langle\cdot,\cdot\rangle$ and $\|\cdot\|$
the standard $L^2$-inner product and $L^2$-norm. For an $S^1$-space $X$, we denote
by $X_{S^1}$ the homotopy quotient of $X$ module the $S^1$-action, i.e.,
$X_{S^1}=S^\infty\times_{S^1}X$. We define the functions
\be \left\{\matrix{[a]=\max\{k\in\Z\,|\,k\le a\}, &
E(a)=\min\{k\in\Z\,|\,k\ge a\} , \cr
                   \varphi(a)=E(a)-[a],   \cr}\right. \lb{1.6}\ee
Specially, $\varphi(a)=0$ if $ a\in\Z\,$, and $\varphi(a)=1$ if $a\notin\Z\,$.
In this article we use only $\Q$-coefficients for all homological modules.

\setcounter{equation}{0}%\setcounter{figure}{0}
\section{Critical point theory for closed characteristics }%{Section 3}

In this section, we describe the critical point theory for closed
characteristics.

As in P.199 of \cite{Eke3}, choose some $\alpha\in(1,\, 2)$ and associate with $U$
a convex function $H_\alpha$ such that $H_\alpha(\lambda x)=\lambda^\alpha H_\alpha(x)$ for $\lambda\ge 0$.
Consider the fixed period problem
\be \left\{\matrix{\dot{x}(t)=JH_\alpha^\prime(x(t)), \cr
     x(1)=x(0).         \cr }\right. \lb{2.1}\ee

Define
\be L_0^{\frac{\alpha}{\alpha-1}}(S^1,\R^{2n})
  =\{u\in L^{\frac{\alpha}{\alpha-1}}(S^1,\R^{2n})\,|\,\int_0^1udt=0\}. \lb{2.2}\ee
The corresponding Clarke-Ekeland dual action functional is defined by
\be \Phi(u)=\int_0^1\left(\frac{1}{2}Ju\cdot Mu+H_\alpha^{\ast}(-Ju)\right)dt,
    \qquad \forall\;u\in L_0^{\frac{\alpha}{\alpha-1}}(S^1,\R^{2n}), \lb{2.3}\ee
where $Mu$ is defined by $\frac{d}{dt}Mu(t)=u(t)$ and $\int_0^1Mu(t)dt=0$,
$H_\alpha^\ast$ is the Fenchel transform of $H_\alpha$ defined by
$H_\alpha^\ast(y)=\sup\{x\cdot y-H_\alpha(x)\;|\; x\in \R^{2n}\}$.
By Theorem 5.2.8 of \cite{Eke3}, $\Phi$ is $C^1$ on
$L_0^{\frac{\alpha}{\alpha-1}}$ and satisfies
the Palais-Smale condition.  Suppose
$x$ is a solution of (\ref{2.1}). Then $u=\dot{x}$ is a critical point
of $\Phi$. Conversely, suppose $u$ is a critical point of $\Phi$.
Then there exists a unique $\xi\in\R^{2n}$ such that $Mu-\xi$ is a
solution of (\ref{2.1}). In particular, solutions of (\ref{2.1}) are in
one to one correspondence with critical points of $\Phi$. Moreover,
$\Phi(u)<0$ for every critical point $u\not= 0$ of $\Phi$.

Suppose $u$ is a nonzero critical point of $\Phi$. Then
the formal Hessian of $\Phi$ at $u$ on $L_0^2(S^1,\R^{2n})$ is defined by
$$ Q(v,\; v)=\int_0^1 (Jv\cdot Mv+(H_\alpha^\ast)^{\prime\prime}(-Ju)Jv\cdot Jv)dt, $$
which defines an orthogonal splitting $L_0^2(S^1,\R^{2n})=E_-\oplus E_0\oplus E_+$ of
$L_0^2(S^1,\; \R^{2n})$ into negative, zero and positive subspaces. The
index of $u$ is defined by $i(u)=\dim E_-$ and the nullity of $u$ is
defined by $\nu(u)=\dim E_0$. Specially $1\le \nu(u)\le 2n$
always holds, cf. P.219 of  \cite{Eke3}.

We have a natural $S^1$-action on $L_0^{\frac{\alpha}{\alpha-1}}(S^1,\; \R^{2n})$ defined by
$\th\cdot u(t)=u(\th+t)$ for all $\th\in S^1$ and $t\in\R$. Clearly
$\Phi$ is $S^1$-invariant. Hence if $u$ is a critical
point of $\Phi$, then the whole orbit $S^1\cdot u$ is formed by
critical points of $\Phi$. Denote by $crit(\Phi)$ the set of
critical points of $\Phi$.
Then $crit(\Phi)$ is compact by the Palais-Smale condition.

For a closed characteristic $(\tau,y)$ on $\Sigma$, we denote by
$y^m\equiv (m\tau, y)$ the $m$-th iteration of $y$ for $m\in\N$.
Let $u^m$ be the unique critical point of $\Phi$ corresponding to $(m\tau, y)$.
Then we define the index $i(y^m)$ and nullity $\nu(y^m)$
of $(m\tau,y)$ for $m\in\N$ by
\be i(y^m)=i(u^m), \qquad \nu(y^m)=\nu(u^m). \lb{2.4}\ee
The mean index of $(\tau,y)$ is defined by
\be \hat{i}(y)=\lim_{m\rightarrow\infty}\frac{i(y^m)}{m}. \lb{2.5}\ee
Note that $\hat{i}(y)>2$ always holds which was proved by Ekeland and
Hofer in \cite{EkH1} of 1987 (cf. Corollary 8.3.2 and Lemma 15.3.2
of \cite{Lon5} for a different proof).

Recall that for a principal $U(1)$-bundle $E\to B$, the Fadell-Rabinowitz index
(cf. \cite{FaR1}) of $E$ is defined to be $\sup\{k\;|\, c_1(E)^{k-1}\not= 0\}$,
where $c_1(E)\in H^2(B,\Q)$ is the first rational Chern class. For a $U(1)$-space,
i.e., a topological space $X$ with a $U(1)$-action, the Fadell-Rabinowitz index is
defined to be the index of the bundle $X\times S^{\infty}\to X\times_{U(1)}S^{\infty}$,
where $S^{\infty}\to CP^{\infty}$ is the universal $U(1)$-bundle.

For any $\kappa\in\R$, we denote by
\be \Phi^{\kappa-}=\{u\in L_0^{\frac{\alpha}{\alpha-1}}(S^1,\R^{2n})\;|\;
             \Phi(u)<\kappa\}. \lb{2.6}\ee
Then as in P.218 of \cite{Eke3}, we define
\be c_i=\inf\{\delta\in\R\;|\: \hat I(\Phi^{\delta-})\ge i\},\lb{2.7}\ee
where $\hat I$ is the Fadell-Rabinowitz index defined above. Then by Proposition 5.3.3
in P.218 of \cite{Eke3}, we have

{\bf Proposition 2.1.} {\it Each $c_i$ is a critical value of $\Phi$. If
$c_i=c_j$ for some $i<j$, then there are infinitely many geometrically
distinct closed characteristics on $\Sg$.}

By Theorem 5.3.4 in P.219 of \cite{Eke3}, we have

{\bf Proposition 2.2.} {\it For every $i\in\N$, there exists a point
$u\in L_0^{\frac{\alpha}{\alpha-1}}(S^1,\R^{2n})$ such that}
\bea
&& \Phi^\prime(u)=0,\quad \Phi(u)=c_i, \lb{2.8}\\
&& i(u)\le 2(i-1)\le i(u)+\nu(u)-1. \lb{2.9}\eea

\setcounter{equation}{0}%\setcounter{figure}{0}
\section{ Index theory for closed characteristics}%{Section 2}

In this section, we recall briefly an index theory for symplectic paths
developed by Y. Long and his coworkers.
All the details can be found in \cite{Lon5}.
These results will be used in the next section.

As usual, the symplectic group $\Sp(2n)$ is defined by
$$ \Sp(2n) = \{M\in {\rm GL}(2n,\R)\,|\,M^TJM=J\}, $$
whose topology is induced from that of $\R^{4n^2}$. For $\tau>0$ we are interested
in paths in $\Sp(2n)$:
$$ \P_{\tau}(2n) = \{\ga\in C([0,\tau],\Sp(2n))\,|\,\ga(0)=I_{2n}\}, $$
which is equipped with the topology induced from that of $\Sp(2n)$. The
following real function was introduced in \cite{Lon3}:
$$ D_{\om}(M) = (-1)^{n-1}\ol{\om}^n\det(M-\om I_{2n}), \qquad
          \forall \om\in\U,\, M\in\Sp(2n). $$
Thus for any $\om\in\U$ the following codimension $1$ hypersurface in $\Sp(2n)$ is
defined in \cite{Lon3}:
$$ \Sp(2n)_{\om}^0 = \{M\in\Sp(2n)\,|\, D_{\om}(M)=0\}.  $$
For any $M\in \Sp(2n)_{\om}^0$, we define a co-orientation of $\Sp(2n)_{\om}^0$
at $M$ by the positive direction $\frac{d}{dt}Me^{t\ep J}|_{t=0}$ of
the path $Me^{t\ep J}$ with $0\le t\le 1$ and $\ep>0$ being sufficiently
small. Let
\bea
\Sp(2n)_{\om}^{\ast} &=& \Sp(2n)\bs \Sp(2n)_{\om}^0,   \nn\\
\P_{\tau,\om}^{\ast}(2n) &=&
      \{\ga\in\P_{\tau}(2n)\,|\,\ga(\tau)\in\Sp(2n)_{\om}^{\ast}\}, \nn\\
\P_{\tau,\om}^0(2n) &=& \P_{\tau}(2n)\bs  \P_{\tau,\om}^{\ast}(2n).  \nn\eea
For any two continuous arcs $\xi$ and $\eta:[0,\tau]\to\Sp(2n)$ with
$\xi(\tau)=\eta(0)$, it is defined as usual:
$$ \eta\ast\xi(t) = \left\{\matrix{
            \xi(2t), & \quad {\rm if}\;0\le t\le \tau/2, \cr
            \eta(2t-\tau), & \quad {\rm if}\; \tau/2\le t\le \tau. \cr}\right. $$
Given any two $2m_k\times 2m_k$ matrices of square block form
$M_k=\left(\matrix{A_k&B_k\cr
                                C_k&D_k\cr}\right)$ with $k=1, 2$,
as in \cite{Lon5}, the $\;\dm$-product of $M_1$ and $M_2$ is defined by
the following $2(m_1+m_2)\times 2(m_1+m_2)$ matrix $M_1\dm M_2$:
$$ M_1\dm M_2=\left(\matrix{A_1&  0&B_1&  0\cr
                               0&A_2&  0&B_2\cr
                             C_1&  0&D_1&  0\cr
                               0&C_2&  0&D_2\cr}\right). \nn$$  %\dm=\diamond
Denote by $M^{\dm k}$ the $k$-fold $\dm$-product $M\dm\cdots\dm M$. Note
that the $\dm$-product of any two symplectic matrices is symplectic. For any two
paths $\ga_j\in\P_{\tau}(2n_j)$ with $j=0$ and $1$, let
$\ga_0\dm\ga_1(t)= \ga_0(t)\dm\ga_1(t)$ for all $t\in [0,\tau]$.

A special path $\xi_n\in\P_{\tau}(2n)$ is defined by
\be \xi_n(t) = \left(\matrix{2-\frac{t}{\tau} & 0 \cr
                                             0 &  (2-\frac{t}{\tau})^{-1}\cr}\right)^{\dm n}
        \qquad {\rm for}\;\quad0\le t\le \tau.  \lb{3.1}\ee

{\bf Definition 3.1.} (cf. \cite{Lon3}, \cite{Lon5}) {\it For any $\om\in\U$ and
$M\in \Sp(2n)$, define
\be  \nu_{\om}(M)=\dim_{\C}\ker_{\C}(M - \om I_{2n}).  \lb{3.2}\ee
For any $\tau>0$ and $\ga\in \P_{\tau}(2n)$, define
\be  \nu_{\om}(\ga)= \nu_{\om}(\ga(\tau)).  \lb{3.3}\ee

If $\ga\in\P_{\tau,\om}^{\ast}(2n)$, define
\be i_{\om}(\ga) = [\Sp(2n)_{\om}^0: \ga\ast\xi_n],  \lb{3.4}\ee
where the right hand side of (\ref{3.4}) is the usual homotopy intersection
number, and the orientation of $\ga\ast\xi_n$ is its positive time direction under
homotopy with fixed end points.

If $\ga\in\P_{\tau,\om}^0(2n)$, we let $\mathcal{F}(\ga)$
be the set of all open neighborhoods of $\ga$ in $\P_{\tau}(2n)$, and define
\be i_{\om}(\ga) = \sup_{U\in\mathcal{F}(\ga)}\inf\{i_{\om}(\beta)\,|\,
                       \beta\in U\cap\P_{\tau,\om}^{\ast}(2n)\}.
               \lb{3.5}\ee
Then
$$ (i_{\om}(\ga), \nu_{\om}(\ga)) \in \Z\times \{0,1,\ldots,2n\}, $$
is called the index function of $\ga$ at $\om$. }

Note that when $\om=1$, this index theory was introduced by
C. Conley-E. Zehnder in \cite{CoZ1} for the non-degenerate case with $n\ge 2$,
Y. Long-E. Zehnder in \cite{LZe1} for the non-degenerate case with $n=1$,
and Y. Long in \cite{Lon1} and C. Viterbo in \cite{Vit2} independently for
the degenerate case. The case for general $\om\in\U$ was defined by Y. Long
in \cite{Lon3} in order to study the index iteration theory (cf. \cite{Lon5}
for more details and references).

For any symplectic path $\ga\in\P_{\tau}(2n)$ and $m\in\N$,  we
define its $m$-th iteration $\ga^m:[0,m\tau]\to\Sp(2n)$ by
\be \ga^m(t) = \ga(t-j\tau)\ga(\tau)^j, \qquad
  {\rm for}\quad j\tau\leq t\leq (j+1)\tau,\;j=0,1,\ldots,m-1.
     \lb{3.6}\ee
We still denote the extended path on $[0,+\infty)$ by $\ga$.

{\bf Definition 3.2.} (cf. \cite{Lon3}, \cite{Lon5}) {\it For any $\ga\in\P_{\tau}(2n)$,
we define
\be (i(\ga,m), \nu(\ga,m)) = (i_1(\ga^m), \nu_1(\ga^m)), \qquad \forall m\in\N.
   \lb{3.7}\ee
The mean index $\hat{i}(\ga,m)$ per $m\tau$ for $m\in\N$ is defined by
\be \hat{i}(\ga,m) = \lim_{k\to +\infty}\frac{i(\ga,mk)}{k}. \lb{3.8}\ee
For any $M\in\Sp(2n)$ and $\om\in\U$, the {\it splitting numbers} $S_M^{\pm}(\om)$
of $M$ at $\om$ are defined by
\be S_M^{\pm}(\om)
     = \lim_{\ep\to 0^+}i_{\om\exp(\pm\sqrt{-1}\ep)}(\ga) - i_{\om}(\ga),
   \lb{3.9}\ee
for any path $\ga\in\P_{\tau}(2n)$ satisfying $\ga(\tau)=M$.}

For a given path $\gamma\in {\cal P}_{\tau}(2n)$ we consider to deform
it to a new path $\eta$ in ${\cal P}_{\tau}(2n)$ so that
\begin{equation}
i_1(\gamma^m)=i_1(\eta^m),\quad \nu_1(\gamma^m)=\nu_1(\eta^m), \quad
         \forall m\in {\bf N}, \label{3.10}
\end{equation}
and that $(i_1(\eta^m),\nu_1(\eta^m))$ is easy enough to compute. This
leads to finding homotopies $\delta:[0,1]\times[0,\tau]\to {\rm Sp}(2n)$
starting from $\gamma$ in ${\cal P}_{\tau}(2n)$ and keeping the end
points of the homotopy always stay in a certain suitably chosen maximal
subset of ${\rm Sp}(2n)$ so that (\ref{3.10}) always holds. In fact,  this
set was first discovered in \cite{Lon3} as the path connected component
$\Omega^0(M)$ containing $M=\gamma(\tau)$ of the set
\begin{eqnarray}
  \Omega(M)=\{N\in{\rm Sp}(2n)\,&|&\,\sigma(N)\cap{\bf U}=\sigma(M)\cap{\bf U}\;
{\rm and}\;  \nonumber\\
 &&\qquad \nu_{\lambda}(N)=\nu_{\lambda}(M)\;\forall\,
\lambda\in\sigma(M)\cap{\bf U}\}. \label{3.11}
\end{eqnarray}
Here $\Omega^0(M)$ is called the {\it homotopy component} of $M$ in
${\rm Sp}(2n)$.

In \cite{Lon3}-\cite{Lon5}, the following symplectic matrices were introduced
as {\it basic normal forms}:
\begin{eqnarray}
D(\lambda)=\left(\matrix{\lm & 0\cr
         0  & \lm^{-1}\cr}\right), &\quad& \lm=\pm 2,\lb{3.12}\\
N_1(\lm,b) = \left(\matrix{\lm & b\cr
         0  & \lm\cr}\right), &\quad& \lm=\pm 1, b=\pm1, 0, \lb{3.13}\\
R(\th)=\left(\matrix{\cos\th & -\sin\th\cr
        \sin\th  & \cos\th\cr}\right), &\quad& \th\in (0,\pi)\cup(\pi,2\pi),
                     \lb{3.14}\\
N_2(\om,b)= \left(\matrix{R(\th) & b\cr
              0 & R(\th)\cr}\right), &\quad& \th\in (0,\pi)\cup(\pi,2\pi),
                     \lb{3.15}\end{eqnarray}
where $b=\left(\matrix{b_1 & b_2\cr
               b_3 & b_4\cr}\right)$ with  $b_i\in\R$ such that $(b_2-b_3)\sin\theta>0$, if $ N_2(\omega, b)$ is
trivial; $(b_2-b_3)\sin\theta<0$, if $ N_2(\omega, b)$ is
non-trivial.

Splitting numbers possess the following properties:

{\bf Lemma 3.3.} (cf. \cite{Lon3} and Lemmas 9.1.5 of \cite{Lon5}) {\it Splitting
numbers $S_M^{\pm}(\om)$ are well defined, i.e., they are independent of the choice
of the path $\ga\in\P_\tau(2n)$ satisfying $\ga(\tau)=M$ appeared in (\ref{3.9}).
For $\om\in\U$ and $M\in\Sp(2n)$, splitting numbers $S_N^{\pm}(\om)$ are constant
for all $N\in\Om^0(M)$. }

{\bf Lemma 3.4.} (cf. \cite{Lon3}, Lemma 9.1.5-9.1.6 and List 9.1.12 of \cite{Lon5})
{\it For $M\in\Sp(2n)$ and $\om\in\U$, there hold
\bea
S_M^{+}(\om) &=& S_M^{-}(\overline{\om}), \qquad \forall\om\in\U.  \lb{3.16}\\
S_M^{\pm}(\om) &=& 0, \qquad {\it if}\;\;\om\not\in\sg(M).  \lb{3.17}\\
(S_{N_1(1,a)}^+(1),\,S_{N_1(1,a)}^-(1)) &=& \left\{\matrix{(1,\,1), &\quad {\rm if}\;\; a\ge 0, \cr
(0,\,0) &\quad {\rm if}\;\; a< 0. \cr}\right. \lb{3.18}\\
(S_{N_1(-1,a)}^+(-1),\,S_{N_1(-1,a)}^-(-1)) &=& \left\{\matrix{(1,\,1), &\quad {\rm if}\;\; a\le 0, \cr
(0,\,0) &\quad {\rm if}\;\; a> 0. \cr}\right. \lb{3.19}\\
(S_{R(\theta)}^+(e^{\sqrt{-1}\theta}),\,S_{R(\theta)}^-(e^{\sqrt{-1}\theta}))
&=& (0,\, 1) \qquad {\it if}\;\;\th\in (0,\pi)\cup(\pi,2\pi).\lb{3.20}\\
(S_{N_2(\om, b)}^+(e^{\sqrt{-1}\theta}),\,S_{N_2(\om, b)}^-(e^{\sqrt{-1}\theta}))
&=& (1,\, 1) \qquad {\it if}\;\; (b_2-b_3)\sin\theta<0..\lb{3.21}\\
(S_{N_2(\om, b)}^+(e^{\sqrt{-1}\theta}),\,S_{N_2(\om, b)}^-(e^{\sqrt{-1}\theta}))
&=& (0,\, 0) \qquad {\it if}\;\; (b_2-b_3)\sin\theta>0..\lb{3.22}
\eea

For any $M_i\in\Sp(2n_i)$ with $i=0$ and $1$, there holds }
\be S^{\pm}_{M_0\dm M_1}(\om) = S^{\pm}_{M_0}(\om) + S^{\pm}_{M_1}(\om),
    \qquad \forall\;\om\in\U. \lb{3.23}\ee

We have the following

{\bf Theorem 3.5.} (cf. \cite{Lon4} and Theorem 1.8.10 of \cite{Lon5}) {\it For
any $M\in\Sp(2n)$, there is a path $f:[0,1]\to\Om^0(M)$ such that $f(0)=M$ and
\be f(1) = M_1\dm\cdots\dm M_k,  \lb{3.24}\ee
where each $M_i$ is a basic normal form listed in (\ref{3.12})-(\ref{3.15})
for $1\leq i\leq k$. In particular, we have $e(f(1))\le e(M)$.}

Let $\Sigma\in\H(2n)$. Using notations in \S1,
for any $(\tau,y)\in\J(\Sigma,\alpha)$ and $m\in\N$, we define
its $m$-th iteration $y^m:\R/(m\tau\Z)\to\R^{2n}$ by
\be y^m(t) = y(t-j\tau), \qquad {\rm for}\quad j\tau\leq t\leq (j+1)\tau,
       \quad j=0,1,2,\ldots, m-1. \lb{3.25}\ee
We still denote by $y$ its extension to $[0,+\infty)$.

We define via Definition 3.2 the following
\bea  S^\pm_y(\om) &=& S^\pm_{\ga_y(\tau)}(\om),  \lb{3.26}\\
  (i(y,m), \nu(y,m)) &=& (i(\ga_y,m), \nu(\ga_y,m)),  \lb{3.27}\\
   \hat{i}(y,m) &=& \hat{i}(\ga_y,m),  \lb{3.28}\eea
for all $m\in\N$, where $\ga_y$ is the associated symplectic path of $(\tau,y)$.
Then we have the following

{\bf Theorem 3.6.} (cf. Lemma 1.1 of \cite{LoZ1}, Theorem 15.1.1 of \cite{Lon5}) {\it Suppose
$(\tau,y)\in \J(\Sigma, \alpha)$. Then we have
\be i(y^m)\equiv i(m\tau ,y)=i(y, m)-n,\quad \nu(y^m)\equiv\nu(m\tau, y)=\nu(y, m),
       \qquad \forall m\in\N, \lb{3.29}\ee
where $i(y^m)$ and $\nu(y^m)$ are the index and nullity
defined in \S2. In particular, (\ref{2.5}) and (\ref{3.8})
coincide, thus we simply denote them by $\hat i(y)$.}

\setcounter{equation}{0}%\setcounter{figure}{0}
\section{Proofs of the main theorems }%{Section 3}

In this section we give the proofs of the main theorems.

Suppose $(\tau,\,y)\in\J(\Sigma,\alpha)$.
Then by Lemma 1.3 of \cite{LoZ1} or Lemma 15.2.4 of \cite{Lon5}, there exist $P_y\in \Sp(2n)$ and $M_y\in \Sp(2n-2)$ such
that
\be \ga_y(\tau)=P_y^{-1}(N_1(1,\,1)\dm M_y)P_y,
   \lb{4.1}\ee
here we use notations in \S3.

Since $H_2(\cdot)$ is positive homogeneous of degree-two,
by (\ref{1.5}) we have
\bea |x|^2R^{-2}\le H_2(x)\le |x|^2r^{-2},\qquad \forall x\in\Sg.
\lb{4.2}\eea

Recall that the action of a closed characteristic $(\tau,\, y)$
is defined by (cf. P190 of \cite{Eke3})
\be A(\tau, y)=\frac{1}{2}\int_0^\tau (Jy\cdot \dot y)dt. \lb{4.3}\ee
Note that $A(\tau, y)$ is a geometric quantity depending only on how many
times one runs around the closed characteristic.
In fact, we have $A(\tau, y)=A(\sg, y\circ\phi)$ for any orientation-preserving
diffeomorphism $\phi: \R/\sg\Z\rightarrow\R/\tau\Z$.

Comparing with the theorem of Morse-Schoenberg in the
study of geodesics, we have the following

{\bf Lemma 4.1.} {\it let $\Sigma\subset\R^{2n}$ be a compact convex
hypersurface which is $(r, R)$-pinched.
Suppose $(\tau, y)$ is a closed characteristic  on $\Sg$.
Then we have the following
\bea &&A(\tau, y)>k\pi R^2\Rightarrow i(y)\ge 2nk,\lb{4.4}\\
&&A(\tau, y)<k\pi r^2\Rightarrow i(y)+\nu(y)\le 2n(k-1)-1.\lb{4.5}
\eea
}

{\bf Proof.} By Proposition 1.7.5 of \cite{Eke3},
we have
\bea i(y)=i_{T_\alpha}(x_\alpha),\qquad\forall\alpha\in (1,\, 2]
\lb{4.6}\eea
where $(T_\alpha,\, x_\alpha)$ is a solution of
\be
\left\{\matrix{\dot{x}(t)=JH_\alpha^\prime(x(t)),
             &&\quad H_\alpha(x(t))=1,\qquad \forall t\in\R. \cr
     x(T)=x(0) && \cr }\right. \lb{4.7}\ee
and $i_{T_\alpha}(x_\alpha)$ is defined in \S 1.6 of \cite{Eke3}.

Now consider the following three Hamiltonian systems
\bea
&&\dot y=2JR^{-2}y,\lb{4.8}\\
&&\dot y=2JH_2^{\prime\prime}(x_2(t))y,\lb{4.9}\\
&&\dot y=2Jr^{-2}y,\lb{4.10}
\eea
and the three corresponding quadratic forms on
$$L_0^2([0,\,s],\,\R^{2n})
  =\{u\in L^2([0,\,s],\,\R^{2n})\,|\,\int_0^1udt=0\}$$
\bea && Q_s^R(v,\; v)=\int_0^s \left(Jv\cdot Mv+\frac{R^2}{2}
\|v\|^2\right)dt,
\lb{4.11}\\
&&Q_s(v,\; v)=\int_0^s (Jv\cdot Mv+H_2^{\prime\prime}(x_2(t))^{-1}Jv\cdot Jv)dt,
\lb{4.12}\\
&& Q_s^r(v,\; v)=\int_0^s \left(Jv\cdot Mv+\frac{r^2}{2}
\|v\|^2\right)dt,
\lb{4.13}
\eea
Note that by (\ref{1.5}) we have
$Q_s^R(v,\; v)\ge Q_s(v,\; v)\ge Q_s^r(v,\; v)$.
Thus we have $i_s^R\le i_s\le i_s^r$, where
$i_s^R,\, i_s$ and $i_s^r$ are the indices of
$Q_s^R,\,Q_s$ and $Q_s^r$.

Note that by (21) in P.191 of \cite{Eke3}, we have
$A(\tau,\,y)=T_2$.
Hence we have
\bea i(y)=i_{T_2}(x_2)\ge i_{T_2}^R\ge 2nk,\lb{4.14}
\eea
where the last inequality follows by
$T_2=A(\tau, y)>k\pi R^2$ and Lemma 1.4.13 of \cite{Eke3}.

Denote by $L_0^2([0,\,T_2],\,\R^{2n})=E_-\oplus E_0\oplus E_+$
the orthogonal splitting  of
$L_0^2([0,\,T_2],\,\R^{2n})$ into negative, zero and positive subspaces.
Then we have the following observation:
If $V$ is a subspace of $L_0^2([0,\,T_2],\,\R^{2n})$
such that $Q_{T_2}$ is negative semi-definite,
i.e., $\xi\in V$ implies $Q_{T_2}(\xi,\,\xi)\le 0$,
then $\dim V\le \dim E_-+\dim E_0$.
In fact, this is a simple fact of linear algebra:
Let
$$pr_-: L_0^2([0,\,T_2],\,\R^{2n})=E_-\oplus E_0\oplus E_+
\rightarrow E_-$$
be the orthogonal projection. Consider
$pr_-|V: V\rightarrow E_-$. Then $\xi\in\ker pr_-|V$
must belong to $E_0\oplus E_+$. That is, since
$Q_{T_2}(\xi,\,\xi)\le 0$, $\xi\in E_0$. From
$$ \dim V=\dim(\im pr_-|V)+\dim(\ker pr_-|V)$$
we prove our claim.

Let $\epsilon>0$ be small enough such that
$A(\tau, y)<k\pi (r-\epsilon)^2$.
If $V$ is a subspace of $L_0^2([0,\,T_2],\,\R^{2n})$
such that $Q_{T_2}|V\le 0$, then
$Q^{r-\epsilon}_{T_2}|V<0$. Hence we have
$\dim V\le i_{T_2}^{r-\epsilon}$.
In particular, we have
$\dim E_-+\dim E_0\le i_{T_2}^{r-\epsilon}$.
Hence by Lemma 1.4.13 of \cite{Eke3},
we have
\bea \dim E_-+\dim E_0\le i_{T_2}^{r-\epsilon}\le 2n(k-1).
\lb{4.15}
\eea
Note that $\dim E_-=i(y)$ and $\dim E_0=\nu(y)+1$.
Hence the lemma follows.\hfill\hb

Suppose $M\in Sp(2n)$ has the normal form decomposition
\bea M=&&N_1(1,1)^{\diamond p_-} \diamond I_{2p_0}\diamond
N_1(1,-1)^{\diamond p_+}
\diamond N_1(-1,1)^{\diamond q_-} \diamond (-I_{2q_0})\diamond
N_1(-1,-1)^{\diamond q_+}\nn\\
&&\diamond R(\theta_1)\diamond\cdots\diamond R(\theta_r)
\diamond N_2(\omega_1, u_1)\diamond\cdots\diamond N_2(\omega_{r_*}, u_{r_*}) \nn\\
&&\diamond N_2(\lm_1, v_1)\diamond\cdots\diamond N_2(\lm_{r_0}, v_{r_0})
\diamond M_0 \lb{4.16}\eea
where $ N_2(\omega_j, u_j) $s are
non-trivial and   $ N_2(\lm_j, v_j)$s  are trivial basic normal
forms; $\sigma (M_0)\cap U=\emptyset$; $p_-$, $p_0$, $p_+$, $q_-$,
$q_0$, $q_+$, $r$, $r_*$ and $r_0$ are non-negative integers;
$\omega_j=e^{\sqrt{-1}\alpha_j}$, $
\lambda_j=e^{\sqrt{-1}\beta_j}$; $\theta_j$, $\alpha_j$, $\beta_j$
$\in (0, \pi)\cup (\pi, 2\pi)$; these integers and real numbers
are uniquely determined by $M$.

We have the following lemma concerning the iteration of
indices.

{\bf Lemma 4.2.} {\it Suppose
$(\tau,y)\in \J(\Sigma, \alpha)$ such that
$\ga_y(\tau)$ can be deformed in $\Omega^0(\ga_y(\tau))$
to $M$ as in (\ref{4.16}).
Then we have $i(y, 2)-2i(y, 1)\le n$ and
$i(y, 2)+\nu(y, 2)-2(i(y, 1)+\nu(y, 1))\ge 1-n$.
In particular, we have the following

(i) if $i(y, 2)-2i(y, 1)=n$, then we have
\bea M=N_1(1,1)^{\diamond p_-} \diamond I_{2p_0}\diamond
 R(\theta_1)\diamond\cdots\diamond R(\theta_r)
\lb{4.17}\eea
with $p_-+p_0+r=n$ and $\theta_k\in (\pi, 2\pi)$
for $1\le k\le r$. In particular,
$(\tau, y)$ is strictly elliptic.

(ii) if $i(y, 2)+\nu(y, 2)-2(i(y, 1)+\nu(y, 1))= 1-n$,
then we have
\bea M=N_1(1,1) \diamond I_{2p_0}\diamond
N_1(1,-1)^{\diamond p_+}\diamond
 R(\theta_1)\diamond\cdots\diamond R(\theta_r)
\lb{4.18}\eea
with $p_0+p_++r=n-1$ and $\theta_k\in (0, \pi)$
for $1\le k\le r$. In particular,
$(\tau, y)$ is strictly elliptic.
}

{\bf Proof.} By the Bott-type formulae, cf.
Theorem 9.2.1 of \cite{Lon5}, we have
\bea &&i(y, 2)=i(\ga_y, 2)=i_1(\ga_y)+i_{-1}(\ga_y)=i(y, 1)+i_{-1}(\ga_y),\lb{4.19}\\
&&\nu(y, 2)=\nu(\ga_y, 2)=\nu_1(\ga_y)+\nu_{-1}(\ga_y)=\nu(y, 1)+\nu_{-1}(\ga_y),\lb{4.20}
\eea
Hence we have $i(y, 2)-2i(y, 1)=i_{-1}(\ga_y)-i_1(\ga_y)$.
By (\ref{3.9}) we have
\bea i_{-1}(\ga_y)-i_1(\ga_y)=S_M^+(1)
+\sum_{0<\theta<\pi}(S_M^+(e^{\sqrt{-1}\theta})-S_M^-(e^{\sqrt{-1}\theta}))
-S_M^-(-1).\lb{4.21}
\eea
Thus $i_{-1}(\ga_y)-i_1(\ga_y)\le n$ and (i) holds
by Lemma 3.4.

Note that we have
\bea \nu(y, 2)-2\nu(y, 1)=\nu_{-1}(\ga_y)-\nu_1(\ga_y).
\lb{4.22}\eea
Thus $i(y, 2)+\nu(y, 2)-2(i(y, 1)+\nu(y, 1))\ge 1-n$
and (ii) holds by (\ref{4.1}) and Lemma 3.4. \hfill\hb

We have the following theorem due to
I. Ekeland and  J. Lasry.

{\bf Theorem 4.3.} {\it Let $U\subset \R^{2n}$ be a
convex compact set with non-empty interior,
and let $\Sg$ be its boundary. Assume there is a point
$x_0\in\R^{2n}$ such that
\bea r\le |x-x_0|\le R,\qquad \forall x\in\Sg
\lb{4.23}\eea
and $\frac{R}{r}<\sqrt{2}$.
Then $\Sg$ carries at least $n$ geometrically distinct
closed characteristics $\{(\tau_1, y_1),\ldots,(\tau_n, y_n)\}$
where $\tau_i$ is the minimal period of $y_i$,
and the actions $A(\tau_i, y_i)$ satisfy:
\bea \pi r^2\le A(\tau_i, y_i)\le \pi R^2,\qquad 1\le i\le n.
\lb{4.24}\eea
}

By the proof of the  above theorem
and Proposition 2.2, we have

{\bf Lemma 4.4.} {\it Assume $\{(\tau_1, y_1),\ldots,(\tau_n, y_n)\}$
are the closed characteristics found in Theorem 4.3.
Then we have
\bea
&& \Phi^\prime(u_i)=0,\quad \Phi(u_i)=c_i, \lb{4.25}\\
&& i(u_i)\le 2(i-1)\le i(u_i)+\nu(u_i)-1, \lb{4.26}\eea
for $1\le i\le n$, where $u_i$ is the unique critical point of
$\Phi$ corresponding to $(\tau_i, y_i)$. \hfill\hb
}

Now we give the proofs of the main theorems.

{\bf Proof of Theorem 1.1.}
Suppose $\Sigma\subset\R^{2n}$ is a compact convex
hypersurface which is $(r, R)$-pinched with
$\frac{R}{r}<\sqrt{\frac{3}{2}}$.
Then we have (\ref{1.5}) and (\ref{4.2}).
From (\ref{4.2}) we have
\bea r\le |x|\le R,\qquad \forall x\in\Sg
\lb{4.27}\eea
Thus by Theorem 4.3, we obtain $n$
geometrically distinct prime
closed characteristics
$\{(\tau_1, y_1),\ldots,(\tau_n, y_n)\}$
such that (\ref{4.24})-(\ref{4.26}) hold.

{\bf Claim.}  {\it The closed characteristics
$(\tau_1, y_1)$ and  $(\tau_n, y_n)$ must be strictly elliptic. }

Note that $i(y_1)=0$ by (\ref{4.26}).
Thus by Theorem 5.1.10 of \cite{Eke3}, we have
$(\tau_1, y_1)$  must be strictly elliptic.
Here we can give another proof.
By (\ref{4.24}) and $A(2\tau_i, y_i)=2A(\tau_i, y_i)$,
we have
\bea 2\pi r^2\le A(2\tau_i, y_i)\le 2\pi R^2,\qquad 1\le i\le n.
\lb{4.28}\eea
Since $\frac{R}{r}<\sqrt{\frac{3}{2}}$,
we have $A(2\tau_i, y_i)\ge 2\pi r^2>\frac{4}{3}\pi R^2$.
Thus by Lemma 4.1, we have
\be i(y_i^2)\ge 2n.\lb{4.29}\ee
Hence by Theorem 3.6 we have
\bea i(y_1, 2)-2i(y_1, 1)=i(y_1^2)+n-2(i(y_1)+n)\ge n.
\lb{4.30}
\eea
Thus by (ii) of Lemma 4.2, we have $(\tau_1, y_1)$ is strictly elliptic.

Note that
\bea i(y_n)\le 2(n-1)\le i(y_n)+\nu(y_n)-1
\lb{4.31}\eea
by (\ref{4.26}).
On the other hand, we have
$A(2\tau_i, y_i)\le 2\pi R^2<3\pi r^2$.
Thus by Lemma 4.1, we have
\be i(y_i^2)+\nu(y_i^2)\le 4n-1.\lb{4.32}\ee
Hence by Theorem 3.6 we have
\bea &&i(y_n, 2)+\nu(y_n, 2)-2(i(y_n, 1)+\nu(y_n, 1))\nn\\
=&&i(y_n^2)+n-2(i(y_n)+n)+\nu(y_n^2)-2\nu(y_n)\le 1-n.
\lb{4.33}
\eea
Thus by (i) of Lemma 4.2, we have $(\tau_n, y_n)$ is strictly elliptic.
\hfill\hb

{\bf Proof of Theorem 1.2.}
Suppose $(\tau ,y)$ is a hyperbolic closed characteristic.
Then we have
\be \ga_y(\tau)=P_y^{-1}(N_1(1,\,1)\dm M_y)P_y,
   \lb{4.34}\ee
with $\sg(M_y)\cap \U=\emptyset$.
Thus by Theorem 8.3.1 of \cite{Lon5} and Theorem 3.6, we have
\bea i(y^m)=m(i(y)+n+1)-n-1,\quad \nu(y^m)=1,\qquad \forall m\in\N.
\lb{4.35}\eea
Now suppose $\{(\tau_1, y_1),\ldots,(\tau_n, y_n)\}$
are the $n$ geometrically distinct prime
closed characteristics obtained in Theorem 4.3.
Thus if $(\tau_i, y_i)$ is hyperbolic, we have
$i(y_i)= 2(i-1)$ by (\ref{4.26})
\bea i(y_i^m)=m(i(y_i)+n+1)-n-1= 2m(i-1)+(m-1)(n+1),\qquad \forall m\in\N.
\lb{4.36}\eea
Hence by (\ref{4.29}) and (\ref{4.32}),  we have
\bea &&2n\le i(y_i^2)=4(i-1)+n+1,
\lb{4.37}\\
&&4(i-1)+n+1+1=i(y_i^2)+\nu(y_i^2)\le 4n-1.\lb{4.38}
\eea
Hence we have
\bea n-1\le 4(i-1)\le 3(n-1).\lb{4.39}
\lb{4.39}\eea
Thus we have
\bea E\left(\frac{n-1}{4}\right)\le (i-1)\le \left[\frac{3(n-1)}{4}\right].\lb{4.39}
\lb{4.40}\eea
Hence there are at most
$[\frac{3(n-1)}{4}]-E\left(\frac{n-1}{4}\right)+1$
hyperbolic closed characteristics in
$\{(\tau_1, y_1),\ldots,(\tau_n, y_n)\}$.
This implies that there are at least
$$n-\left[\frac{3(n-1)}{4}\right]+E\left(\frac{n-1}{4}\right)-1
=2\left[\frac{n+2}{4}\right]$$
non-hyperbolic closed characteristics
on $\Sg$. \hfill\hb

\bibliographystyle{abbrv}

\medskip

\end{document}